\documentclass[11pt]{article}
\usepackage{amsmath,amsthm,amssymb,hyperref}

\newcommand{\al}{\alpha}

\newcommand{\ep}{\epsilon}

\newcommand{\bR}{\mathbb{R}}

\newcommand{\bN}{\mathbb{N}}

\newcommand{\pt}{\partial_t}
\newcommand{\pa}{\partial}

\newcommand{\beq}{\begin{equation}}
\newcommand{\eneq}{\end{equation}}

\newtheorem{thm}{Theorem}

\newtheorem{prop}[thm]{Proposition}
\newtheorem{rem}{Remark}

\numberwithin{equation}{section}

\begin{document}

\title{Sharp Global Existence for Semilinear Wave Equation with Small
Data\thanks{Supported by NSF of China 10571158 and Zhejiang
Provincial NSF of China (Y605076)}}
\author{Daoyuan Fang and Chengbo Wang\thanks{email: DF:
dyf@zju.edu.cn, CW: wangcbo@yahoo.com.cn}\\
Department of Mathematics, Zhejiang University,\\ Hangzhou, 310027,
China}

\maketitle


 The global existence in time for nonlinear wave equation with small
data usually require high Sobolev regularity, when one dealt with
them by classical energy method (see \cite{Ch85}, \cite{Kl85} for
example). The purpose of this note is to give the sharp regularity
global existence for semilinear equation with the power nonlinearity
of the derivative, the counterpart of quasilinear equation or the
quadratic nonlinearity seems still unreachable.

Consider the following Cauchy problem(denote $\Box:=\pt^2-\Delta$
and $\pa=(\pt,\pa_x)$) \beq\label{SLW}\left\{\begin{array}{l} \Box u
= \sum_{|\al|=k} c_\al (\pa u)^\al := N(u) \\ u(0,x)=u_0\in H^s, \
\pt u(0,x) = u_1 \in H^{s-1} \end{array}\right. \eneq Let
$s_c=\frac{n+2}{2}-\frac{1}{k-1}$ be the scaling index, we have
\begin{thm}\label{fw5-thm-Glob}
Let $\|u_0\|_{H^s}+\|u_1\|_{H^{s-1}}\le \ep$ with $\ep$ small
enough, and \beq\label{sreq}\left\{\begin{array}{ll}s>s_c& {\rm if}\
k-1=\frac{4}{n-1}\vee 2\ {\rm and}\ n\neq 3\\
s\ge s_c & {\rm if}\  k-1>\frac{4}{n-1}\vee 2,
\end{array}\right.\eneq then the equation \eqref{SLW} has a
unique global solution in $C_t H^s$ such that $\pa u\in L_t^\infty
H^{s-1}\cap L_t^{k-1}L^\infty$. Moreover, if $k=n=3$, then the
lifespan $T_*$ of the solution with $s>2$ is at least of order
$\exp(c \ep^{-2})$ with $c\ll 1$.
\end{thm}

We will prove a similar result for the initial data which are
spherical symmetric in addition. For such purpose, we introduce a
concept here. We say that the equation \eqref{SLW} is {\bf radial},
if $u(t,x)$ is any solution of the equation, then for any rotation
$S$ in $\bR^n$, $u(t,Sx)$ is still a solution of the same equation.
For example, when $k=2$, the radial equation must take the form of
$$\Box u=c_1 (\pt u)^2+ c_2  |\nabla u|^2\ .$$

\begin{thm}\label{fw5-thm-RadialGlob}
Let $n\ge 2$ and $k> \frac{n+1}{n-1}\vee 2$, and consider the radial
equation, then there exists a global solution in time for $s\ge s_c$
with small radial data.
\end{thm}

\begin{rem}
  The requirement for regularity in Theorem \ref{fw5-thm-Glob} and
  \ref{fw5-thm-RadialGlob} are essentially sharp. Since for the
  equation
$$\Box  u= |\pt u|^{k-1} \pt u\ ,$$
it's well-known that the problem is ill posed in $H^s$ for $s<s_c$
(see Theorem 2 in \cite{FW3} for example), in the sense that,  there
is a sequence of data $f_j, g_j\in C^\infty_0 (B_{R_j})$, for which
the lifespan of the solutions $u_j$ tends to zero as the data's norm
and $R_j$ goes to $0$, under the condition that the solutions obey
finite speed of propagation. Note that the initial data $f_j, g_j$
can be radial functions. Thus for such $s$, we can not hope any
existence results as in these Theorems.
\end{rem}

\begin{rem} For the case $n=k=3$, we have almost global existence in
general and global existence for the radial data. Thus a natural
question is: To what extent does the result of global existence
depend on the radial symmetry? The answer is that it is very little.
In fact, in \cite{MaNaNaOz05}, the authors show that for any small
data with additional rotation regularity, there is global existence
for the equation \eqref{SLW}.
\end{rem}

\begin{rem}
It's regret that such argument can not apply to the more interesting
case $k=2$, since it's well known that the corresponding $L^1
L^\infty$ Strichartz estimate is not hold true in general. For the
local result for semilinear and quasilinear equation, one can refer
to \cite{Ta99}, \cite{SmTa05} and references therein.
\end{rem}

We will use the Strichartz estimates to prove the result. For the
details of the Strichartz estimates, one may consult \cite{FW2} and
references therein.
\begin{prop}[Strichartz Estimate]\label{stri}
Let $u$ be the solution of the linear wave equation and $q< \infty$, 
then for $(q,n)\neq (2,3)$ \beq\label{StrichartzEstimate}\|\pa
u\|_{L^q L^{\infty}\cap L^{\infty}H^{s-1}}\le C_q \|\pa
u(0)\|_{H^{s-1}}\eneq with $s\ge \frac{n+2}{2}-\frac{1}{q}$ and $q>
\frac{4}{n-1}\vee 2$ or $s> \frac{n+2}{2}-\frac{1}{q}$ and $q=
\frac{4}{n-1}\vee 2$. For the case $(q,n)=(2,3)$ and $s>2$, we have
\beq\label{endptSEst}\|\pa u\|_{L^2([0,T],
L^{\infty})}+(\ln(1+T))^{1/2} \|\pa
u\|_{L^{\infty}([0,T],H^{s-1})}\le C (\ln(1+T))^{1/2} \|\pa
u(0)\|_{H^{s-1}}.\eneq Moreover, if $u$ is spatial radial function,
then we have \eqref{StrichartzEstimate} with $s\ge
\frac{n+2}{2}-\frac{1}{q}$ for all $q>\frac{2}{n-1}$ and $q\ge 2.$
\end{prop}

We'll use Picard's iteration argument to give the proof. First, we
give the proof for the case $k-1\ge \frac{4}{n-1}\vee 2$ and
$(n,k)\neq (3,3)$.

Let $u^{(0)}=0$ and then define $u^{(m+1)}$ ($m\in\bN$) to be the
solution of the problem
$$\Box u^{(m+1)}=N(u^{(m)})$$ with the given data $(u_0, u_1)$.
We'll see below that $(\pt u^{(m)}, \pa_x u^{(m)})$ is a Cauchy
sequence in $C_t H^{s-1}\cap L_t^{k-1}L^\infty$ with  the norm
$L_t^\infty H^{s-1}\cap L_t^{k-1}L^\infty$ if
$\|u_0\|_{H^s}+\|u_1\|_{H^{s-1}}=\ep$ is small enough.

We claim  that for any $m\in \bN$, $u^{(m)}\in C H^s\cap C^1
H^{s-1}$ and \beq\label{inBall}\|\pa u^{(m)}\|_{L^\infty H^{s-1}\cap
L^{k-1}L^\infty}\le M \ep\eneq with $M$ large enough. In fact, it's
true for $m=0$, and we assume it's true for some $m$, then by
Proposition \ref{stri} with $q=k-1$ and $s$ as in \eqref{sreq},
$$
\begin{array}{lcl}
\|\pa u^{(m+1)}\|_{L^\infty H^{s-1}\cap L^{k-1}L^\infty} &\le&
C(\ep+\|N(u^{(m)})\|_{L^1 H^{s-1}})
\\
 &\le &
 C(\ep+\|\pa u^{(m)}\|_{L^{k-1}L^\infty}^{k-1} \|\pa u^{(m)}\|_{L^\infty
 H^{s-1}})\\
&\le & C(\ep + (M\ep)^k)\le M \ep.
\end{array}
$$ Thus we get \eqref{inBall} by induction.

Now we show that $(\pt u^{(m)}, \pa_x u^{(m)})$ is a Cauchy
sequence in $C_t H^{s-1}\cap L_t^{k-1}L^\infty$  with norm
$L_t^\infty H^{s-1}\cap L_t^{k-1}L^\infty$. Note that for any
$m\in\bN_+$, $u^{(m+1)}-u^{(m)}$ is the solution of equation
$$\Box (u^{(m+1)}-u^{(m)})=N(u^{(m)})-N(u^{(m-1)})$$ with the null
data. Then
$$\begin{array}{lcl}\|\pa(u^{(m+1)}-u^{(m)})\|_{L^\infty H^{s-1}\cap
L^{k-1}L^\infty}&\le& C \|N(u^{(m)})-N(u^{(m-1)})\|_{L^1 H^{s-1}}\\
&\le& C \ep^{k-1} \|\pa(u^{(m)}-u^{(m-1)})\|_{L^\infty H^{s-1}\cap
L^{k-1}L^\infty}\\
&\le& \frac{1}{2}\|\pa(u^{(m)}-u^{(m-1)})\|_{L^\infty H^{s-1}\cap
L^{k-1}L^\infty}.
\end{array}
$$
Thus we have
$$\|\pa(u^{(m+1)}-u^{(m)})\|_{L^\infty H^{s-1}\cap
L^{k-1}L^\infty}\le 2^{-m} \|\pa(u^{(1)}-u^{(0)})\|_{L^\infty
H^{s-1}\cap L^{k-1}L^\infty}\le 2^{-m} M \ep$$ by induction and
\eqref{inBall}. So
\beq\label{Cauchy}\|\pa(u^{(m)}-u^{(l)})\|_{L^\infty H^{s-1}\cap
L^{k-1}L^\infty}\le 2^{1-\max(m,l)} M \ep.\eneq

Therefore, there exist $u^i$, $i\in \{0,1,\cdots,n\}$, such that
$$\pa_{i} u^{(m)}\rightarrow u^i\ \mathrm{in}\ C H^{s-1}\cap
L^{k-1}L^\infty\ .$$ Now we define $$u(t)=u_0+\int_0^t u^0\in
CH^{s-1}.$$ Since
$$u^{(m)}(t)=u_0+\int_0^t \pt u^{(m)},$$
thus for any $0 < T<\infty$, $t\in [0,T],$
$$\pa_i u^{(m)}(t)=\pa_i u_0+\int_0^t \pa_i \pt u^{(m)}\rightarrow
\pa_i u_0+\int_0^t \pa_i u^0=\pa_i u(t)\ \mathrm{in}\ C([0,T],
H^{s-2})$$ and so $\pa_i u=u^i$,
$$\pa_{i}
u^{(m)}\rightarrow \pa_i u\ \mathrm{in}\ C H^{s-1}\cap
L^{k-1}L^\infty\ .$$ Then we can get the solution $u\in CH^s\cap
C^1 H^{s-1}$ of equation \eqref{SLW}.

For the uniqueness and continuous dependence of the initial data,
it's essentially as the above proof. Let $\|(u_0,u_1)\|_{H^s\times
H^{s-1}}\le \ep$ and $\|(v_0,v_1)\|_{H^s\times H^{s-1}}\le \ep$.
Assume $u$ and $v$ are two solutions of \eqref{SLW} with data
$(u_0,u_1)$ and $(v_0,v_1)$ respectively, then $u-v$ is the
solution of equation
$$\Box (u-v)=N(u)-N(v)$$ with the data $(u_0-v_0,u_1-v_1)$.
$$\begin{array}{lcl}\|\pa(u-v)\|_{L^\infty H^{s-1}\cap
L^{k-1}L^\infty}&\le& C(\|(u_0-v_0,u_1-v_1)\|_{H^s\times H^{s-1}}  + \|N(u)-N(v)\|_{L^1 H^{s-1}})\\
&\le& C(\|(u_0-v_0,u_1-v_1)\|_{H^s\times H^{s-1}}  + \ep^{k-1}
\|\pa(u-v)\|_{L^\infty H^{s-1}\cap
L^{k-1}L^\infty})\\
&\le&C\|(u_0-v_0,u_1-v_1)\|_{H^s\times H^{s-1}}+
\frac{1}{2}\|\pa(u-v)\|_{L^\infty H^{s-1}\cap L^{k-1}L^\infty}.
\end{array}
$$
Thus we have \beq\label{unique}\|\pa(u-v)\|_{L^\infty H^{s-1}\cap
L^{k-1}L^\infty}\le C\|(u_0-v_0,u_1-v_1)\|_{H^s\times H^{s-1}}\eneq

This complete the proof for the case $k-1\ge \frac{4}{n-1}\vee 2$
and $(n,k)\neq (3,3)$.

For the case $n=k=3$, it remains to claim alternatively that
\beq\label{induction} \|\pa u^{(m)}\|_{L^\infty_{[0,T]} H^{s-1}}\le
M\ep,\ \|\pa u^{(m)}\|_{L^{2}_{[0,T]}L^\infty}\le c\ll 1\eneq if
$\ln(1+T)\ll \ep^{-2}$. In fact, let $$A_m:=\|\pa
u^{(m)}\|_{L^2([0,T], L^{\infty})}+(\ln(1+T))^{1/2} \|\pa
u^{(m)}\|_{L^{\infty}([0,T],H^{s-1})}\ ,$$ then by inductive
assumption,
$$
\begin{array}{lcl}
A_{m+1}&\le& C \ln(1+T)^{\frac{1}{2}}
(\ep+\|N(u^{(m)})\|_{L^1_{[0,T]} H^{s-1}})
\\
 &\le &
 C\ln(1+T)^{\frac{1}{2}}(\ep+\|\pa u^{(m)}\|_{L^{2}_{[0,T]}L^\infty}^{2} \|\pa u^{(m)}\|_{L^\infty_{[0,T]}
 H^{s-1}})\\
&\le & C \ln(1+T)^{\frac{1}{2}}(\ep + c^2 M\ep)\\
&\le& M \ep \ln(1+T)^{\frac{1}{2}}\ll 1.
\end{array}
$$ Thus we have \eqref{induction} for any $m$.

For the radial cases, it only needs to replace the usual Strichartz
estimate by the required radial $L^{k-1} L^\infty$ estimate  in
Proposition \ref{stri}.


\end{document}